\documentclass[runningheads]{llncs}
\usepackage{graphicx}
\usepackage{booktabs} % For formal tables
\usepackage{amsmath,amssymb, graphicx}
\usepackage[ruled,vlined,linesnumbered]{algorithm2e}
\usepackage{algorithmic}
\SetKwComment{Comment}{$\triangleright$\ }{}
\usepackage[mathscr]{euscript}
\usepackage{color}
\usepackage{url}
\usepackage{mathtools}
\usepackage{enumitem}
\usepackage{hyperref}
\usepackage{xspace}
\usepackage{footnote}
\usepackage{multirow}
\usepackage{times}
\usepackage{bm}
\usepackage{comment}
\usepackage{subcaption}
\DeclareGraphicsExtensions{.pdf,.png,.jpg}

\usepackage[misc]{ifsym}

\newcommand{\hide}[1]{}

\newcommand{\omp}{\textsc{OpenMP}\xspace}
\newcommand{\arma}{\textsc{Armadillo}\xspace}

\newcommand{\T}[1]{\boldsymbol{\mathscr{#1}}}   %% Tensor macro
\newcommand{\tensor}[1]{\boldsymbol{\mathscr{#1}}}   %% Tensor macro
\newcommand{\mat}[1]{\mathbf{#1}}

\newcommand{\argmin}[1]{\underset{#1}{\operatorname{arg}\,\operatorname{min}}\;}

\newcommand{\vest}{\textsc{VeST}\xspace}
\newcommand{\vesta}{\textsc{VeST$_{*}^{auto}$}\xspace}
\newcommand{\vestm}{\textsc{VeST$_{*}^{man}$}\xspace}
\newcommand{\vestf}{\textsc{VeST$_{L_F}^{*}$}\xspace}
\newcommand{\vestl}{\textsc{VeST$_{L_1}^{*}$}\xspace}
\newcommand{\vestaf}{\textsc{VeST$^{auto}_{L_F}$}\xspace}

\newcommand{\vestal}{\textsc{VeST$^{auto}_{L_1}$}\xspace}
\newcommand{\vestml}{\textsc{VeST$^{man}_{L_1}$}\xspace}

\begin{document}
\title{\vest: Very Sparse Tucker Factorization of Large-Scale Tensors }
%
%\title{Contribution Title\thanks{Supported by organization x.}}
\author{MoonJeong Park\inst{1} \and
        Jun-Gi Jang\inst{2} \and
        Lee Sael(\Letter)\inst{2}%\orcidID{0000-0002-9066-5756}
}
%
%\author{First Author\inst{1}\orcidID{0000-1111-2222-3333} \and
%Second Author\inst{2,3}\orcidID{1111-2222-3333-4444} \and
%Third Author\inst{3}\orcidID{2222--3333-4444-5555}}
%
\institute{
Daegu Gyeongbuk Institute of Science and Technology, Daegu, Korea \and
%\email{moonjeong94@dgist.ac.kr} \and
Department of Computer Science and Engineering, Seoul National University, Seoul, Korea
\email{saellee@snu.ac.kr}
%\email{\{elnino4,saellee\}@snu.ac.kr}
}
\authorrunning{M. Park et al.}
\maketitle              % typeset the header of the contribution
\begin{abstract}
% 15~250 words
% The abstract should briefly summarize the contents of the paper in 15--250 words.
% Motivating Question
Given a large tensor, how can we decompose it to sparse core tensor and factor matrices such that it is easier to interpret the results? How can we do this without reducing the accuracy?
Existing approaches either output dense results or give low accuracy.
%
% Solution approach
In this paper, we propose \vest, a tensor factorization method for partially observable data to output a very sparse core tensor and factor matrices.
\vest performs initial decomposition, determines unimportant entries in the decomposition results, removes the unimportant entries, and updates the remaining entries.
To determine unimportant entries, we define and use entry-wise `responsibility' for the decomposed results.
The entries are updated iteratively using a carefully derived coordinate descent rule in parallel for scalable computation.
\vest also includes an auto-search algorithm to give a good tradeoff between sparsity and accuracy.
%
% Result - success of solution approach
Extensive experiments show that our method \vest is at least $2.2$ times sparser and at least $2.8$ times more accurate compared to competitors.
Moreover, \vest is scalable in terms of dimensionality, number of observable entries, and number of threads. 
Thanks to \vest,
we successfully interpret the decomposition result of real-world tensor data based on the sparsity pattern of the factor matrices.
% show the interpretability of the results on real data.

\end{abstract}
\section{Introduction}
\label{sec:intro}
% Definition & importance of the problem
%\red{16 pages Max }
How can we factorize a large tensor to sparse core tensor and factor matrices such that outputs are easier to interpret? How can we do this without sacrificing accuracy?
A tensor is a powerful tool for representing multi-modal data.
Tensor factorization outputs a core tensor and factor matrices which reveal the latent relation of the data.
Tensor factorization can also be viewed as a tool for multi-linear regression problem where only the target values, i.e., values of input tensor entries, are known.
In this perspective, the columns of factor matrices act as latent components, their values as latent feature values, and the cells of the core tensor as weights of the relations between the latent components \cite{kolda2009tensor}.
Sparse tensor factorization aims to output sparse core tensor and factor matrices.
As sparse linear regression model enhances its interpretability \cite{DBLP:conf/kdd/Ribeiro0G16},
sparse factor matrices and a core improve interpretability.

% Challenges of the problem / weaknesses of existing methods
There are two widely used approaches for sparse factorization.
The first approach adds an $L_1$ norm as sparse regularizer to the factorization objective function \cite{morup2008algorithms,madrid2017tensor}.
However, the sparsity is sensitive to lambda values.
The second approach removes elements with small values from the core tensor or the factor matrices~\cite{Sun2017,allen2012sparse,Yi2017}.
However, removing such elements does not necessarily leads to small reconstruction error, and thus value-based pruning sacrifices accuracy.

Other application-specific approaches include utilizing domain-specific knowledge as sparsity constraints~\cite{Lee2018a}, constructing factor matrix from sparse input sampling~\cite{Lee2018b,Mahoney2008}, using smoothing matrices \cite{Pascual-Montano2006,kim2007nonnegative}, and learning sparse dictionary for image data~\cite{Qi2016,Zemin2016,Jiang2018}.
These methods are based on several strong assumptions.
The first approach assumes that prior classification of each mode, e.g., gene sets in omics data, are known.
The second approach assumes that input tensors are already sparse and interpretable, e.g., network data.
The third and fourth approaches assume values in input tensor are smooth, and unsmoothing them does not affect the overall meaning of the data, e.g., image data.
However, these strong assumptions limit their use in general tensor factorizations.

In this paper, we propose \vest(\textsc{Ve}ry \textsc{S}parse \textsc{T}ucker factorization), a scalable and accurate Tucker factorization method to generate sparse factors and a core tensor for large-scale partially observable input tensor.
\vest outputs very sparse factorization results by carefully determining the importance of elements of factors and the core, and pruning unimportant ones. \vest guarantees that the sparsity non-decreases in the update process by carefully derived update rules.
Often, increasing sparsity too much degrades accuracy.
\vest gives an algorithm to automatically determine a reasonable sparsity which offers a good balance with regard to accuracy.
The very sparse result of \vest helps interpreting the result of Tucker factorization and easily revealing the relations of dimensions in multilinear regression.

% Table~\ref{tab:summary} shows the comparison of \vest and standard tensor factorization methods.
% \vest is the only method that satisfies all of the desired properties: sparsity, accuracy, scalability, and memory efficiency.

% % Our Contributions % %
\begin{figure}[t!]
    \centering
     \includegraphics[width=0.8\textwidth]{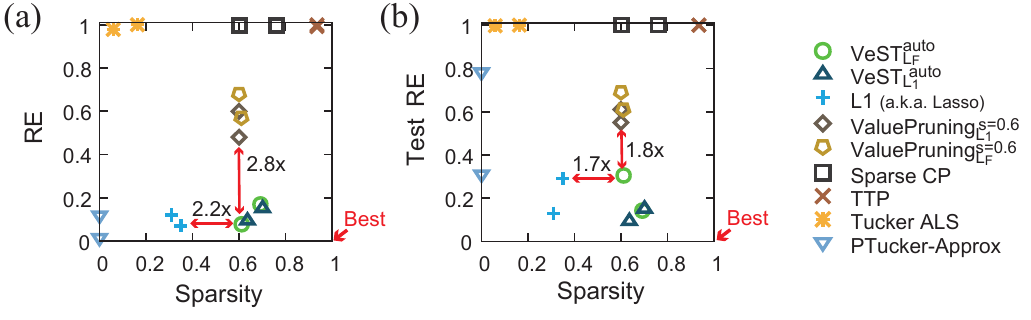}
    \caption{Sparsity and accuracy of \vest and competitors on Yelp-s and AmazonFood-s.
    \vest generates sparse and accurate results that generalize well on unseen data: points are located in the bottom right areas of both RE (reconstruction error) and Test RE plots.
    %    Measurements are averages of five runs.
    }
    \label{fig:accuracy}
\end{figure}

\noindent Our main contributions are as follows:
\begin{itemize}[noitemsep,topsep=0pt]
    \item \textbf{Algorithm.} We propose \vest, an efficient Tucker factorization method for partially observable data, which produces very sparse outputs for better interpretability without loss of accuracy.
        \vest also provides an algorithm to automatically determine a sparsity which gives a good trade-off with regard to accuracy.
    \item \textbf{Theory.} We carefully derive parallelizable coordinate decent update rules for core tensor and factor matrices, and prove their correctness. We also analyze the time and space complexity of \vest.
    \item \textbf{Performance.} \vest provides better sparsity, accuracy, and scalability compared to other methods (see Figure~\ref{fig:accuracy}).
\end{itemize}
%The rest of this paper is organized as follows. Section 2 explains preliminaries and related works.
%Section 3 describes our proposed method \vest.
%Section 4 presents experimental results and discoveries.
%We conclude in Section~\ref{sec:conclusion}.
\noindent The codes and datasets used in this paper are available at \url{http://github.com/leesael/VeST}.

\section{Preliminaries and Related Works}
\label{sec:preliminary}

We introduce concepts of tensor and its operations, % in Section \ref{sec:prelim:tensor}, %tensor operation in Section \ref{sec:prelim:operation},
Tucker factorization, %in Section %\ref{sec:prelim:tf},
and the standard algorithm for Tucker. %in Section \ref{sec:prelim:ALS}.
Table~\ref{tab:Symbols} lists the symbols used.

\begin{table}[htp]
    \centering
    \vspace{-1mm}
    \caption{Table of symbols and definitions.}
    \begin{tabular}{c  p{5cm}| c p{5cm}}
 %       \hline
        %\textbf{Symbol} & \textbf{Definition}  & \textbf{Symbol} & \textbf{Definition} \\
        \hline
        $\tensor{X}$ & input tensor $(\in \mathbb{R}^{I_{1} \times ... \times I_{N}})$ &
        $\tensor{G}$ & core tensor $(\in \mathbb{R}^{J_{1}\times ... \times J_{N}})$ \\
        $N$ & order of $\tensor{X}$ &
        $I_{n},J_n$ & dimensionality of the $n$th mode of $\tensor{X}$ and $\T{G}$, respectively\\
        $\mat{A}^{(n)}$ & $n$th factor matrix $(\in \mathbb{R}^{I_{n} \times J_{n}})$ &
        $a^{(n)}_{i_{n}j_{n}}$ & $(i_{n},j_{n})$th element of $\mat{A}^{(n)}$ \\
        $\Omega$ & set of observable entries of $\tensor{X}$ &
        $|\Omega|$ & number of observable entries of $\tensor{X}$ \\
        $\Omega_{i_n}^{(n)} $ & set of observable entries whose $n$th mode index is $i_n$  &
        $\lambda$ & regularization parameter for core and factor matrices \\
        $\|\tensor{X}\|_F$ & Frobenius norm of tensor $\T{X}$ &
        $\|\tensor{X}\|_1$ & sum of absolute values of tensor $\T{X}$ \\
        $\alpha$ & an entry $(i_1,...,i_N)$ of input tensor $\T{X}$  &
        $\beta$ & an element $(j_1,...,j_N)$ of core tensor  $\T{G}$ \\
        $\alpha_{i_n=i}$ & an entry $(i_1,...,i_n=i,...,i_N)$ of input tensor $\T{X}$  &
        $\beta_{j_n=j}$ & an element $(j_1,...,j_n=j,...,j_N)$ of core tensor  $\T{G}$ \\
%                $T$ & number of CPU cores &
        \hline
    \end{tabular}
    \label{tab:Symbols}
\end{table}

\subsection{Tensor and its Operations} \label{sec:prelim:tensor}
Tensor is multi-dimensional array that contains numbers.
An `order' or `mode' is the number of tensor dimensions, where a $1^{st}$-order tensor represents a vector and a $2^{nd}$-order tensor represents a matrix.
We denote vectors by boldface lowercase letters (e.g., \textbf{a}), matrices by boldface capital letters (e.g., \textbf{A}), and three or higher order tensors by boldface Euler script letters (e.g., \bm{{$\mathscr{X}$}}).
An entry of a $3^{rd}$-order tensor can be expressed with three indices. % (Figure~\ref{fig:3orderTF}).
For example, the ($i_1$, $i_2$, $i_3$)th entry of a $3^{rd}$-order tensor $\T{X} \in \mathbb{R}^{I_{1} \times I_{2} \times I_{3}}$ is denoted by $x_{i_1i_2i_3}$, where index $i_n$ spans from 1 to $I_n$.
%A fiber of a tensor is similar to a column or a row of a matrix.
%A column of matrix is expressed as mode-1 fiber, row as mode-2 fiber. There are three kinds of fiber in 3-order tensor denoted by $x_{:i_2i_3}$, $x_{i_1:i_3}$, and $x_{i_1i_2:}$.

% \subsection{Tensor Operations}
% \label{sec:prelim:operation}
% We review essential tensor operations used in tensor factorization. More tensor operations are summarized in~\cite{kolda2009tensor}.

The size of a tensor is often evaluated by the Frobenius norm. Given an $N$-order tensor $\T{X}$ $(\in \mathbb{R}^{I_{1} \times ... \times I_{N}})$, the \textit{Frobenius norm} of $\T{X}$ is $ ||\T{X}||_F = \sqrt{\sum_{\forall{\alpha}\in\T{X}}{\T{X}^{2}_{\alpha}}}$, where $\alpha\ = (i_1,\cdots, i_N)$ is an index to an entry of input tensor $\T{X}$.
Tensor decomposition often involves matricization of a tensor, and product between a tensor and a matrix.
The \textit{mode-n matricization} of a tensor $\T{X} \in \mathbb{R}^{I_1 \times\cdots \times I_N}$ is denoted as $\mathbf{X}_{(n)}$ and the mapping from an entry $(i_1,\cdots, i_N)$ of $\T{X}$ to an entry $(i_n,j)$ of $\mathbf{X}_{(n)}$ is given by
 $ j=1+\sum_{k=1, k \neq n}^{N}[(i_k-1) \prod_{ m=1 , m \neq n}^{k-1} I_m]$.
 Also, the \textit{n-mode product} of a tensor $\T{X} \in \mathbb{R}^{I_1 \times \cdots \times I_N}$ with a matrix $\mathbf{U} \in \mathbb{R}^{J \times I_n}$ is denoted by $\tensor{X}\times_{n}\mathbf{U}$ ($\in \mathbb{R}^{I_1 \times \cdots \times I_{n-1} \times J \times I_{n+1} \times \cdots \times I_N}$).
Entry-wise, we have
$    {(\tensor{X}\times_{n}\mathbf{U})}_{i_1 \cdots i_{n-1} j i_{n+1} \cdots i_N} =$ $  \sum_{i=1}^{I_n} (\tensor{X}_{\alpha_{i_n=i}} {u}_{ji}) $.
%An n-mode product enables multiplication between a tensor and a matrix.

\subsection{Tucker Factorization}
\label{sec:prelim:tf}
\begin{figure}[h!]
    \centering
     \includegraphics[height=3cm]{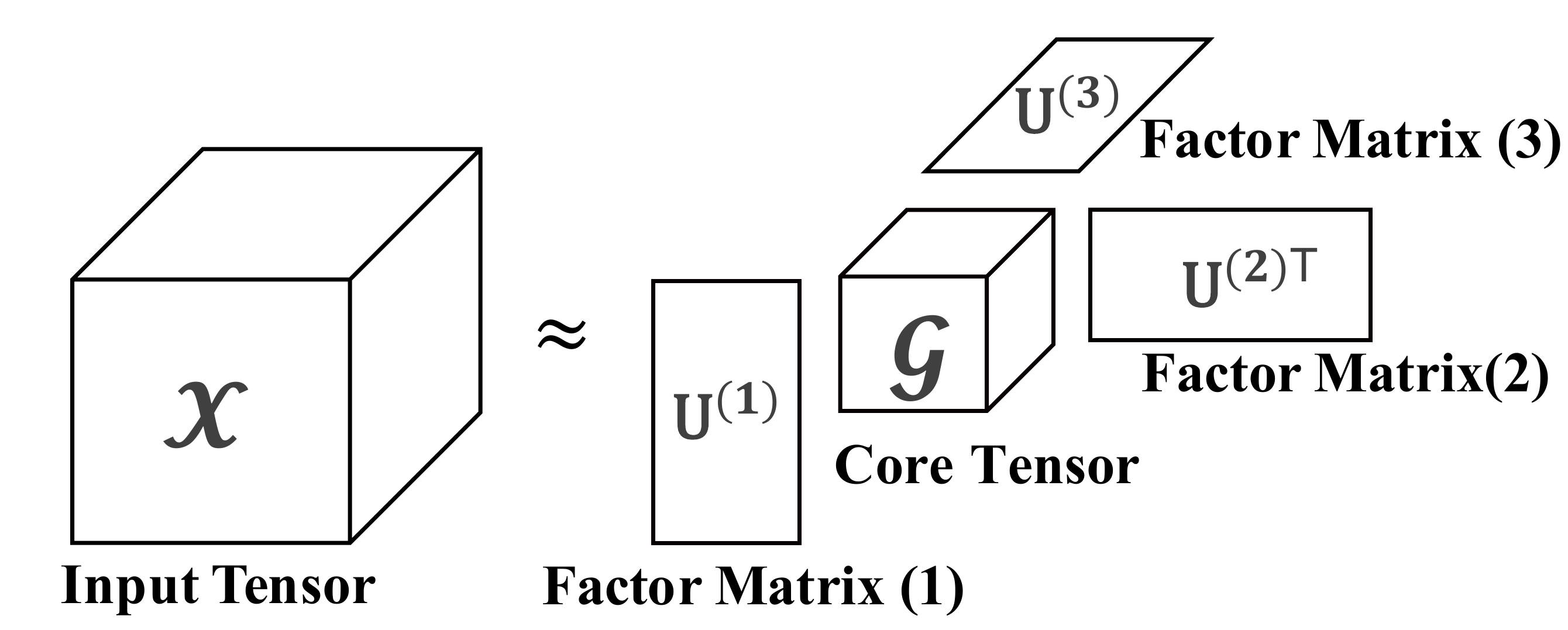}
    \caption{Tucker factorization of a 3-way tensor.}
    \label{fig:3way_tucker}
\end{figure}
Our proposed method \vest is built on top of Tucker factorization, one of the most popular tensor  factorization methods. Given an $N^{th}$-order tensor $\T{X}$ $(\in \mathbb{R}^{I_{1} \times \cdots \times I_{N}})$, Tucker factorization approximates $\T{X}$ by a core tensor $\T{G}$ $(\in \mathbb{R}^{J_{1}\times \cdots \times J_{N}})$ and factor matrices $\{\mathbf{A}^{(n)}\in \mathbb{R}^{I_{n} \times J_n} | n=1\cdots N\}$ by minimizing the full reconstruction error:
  $\min_{\T{G},\mathbf{A}^{(1)},...,\mathbf{A}^{(N)}} \| \T{X} - \T{G}\times_{1}\mathbf{A}^{(1)} \cdots \times_{N}\mathbf{A}^{(N)}\|_F$.

Figure~\ref{fig:3way_tucker} illustrates a Tucker factorization result for a $3^{rd}$-order tensor. Typically, a core tensor $\T{G}$ is assumed to be smaller and denser than the input tensor $\T{X}$. %, and factor matrices $\mathbf{A}^{(n)}$ are orthogonal. %normally orthogonal but non-unique.
%Note that orthogonality alleviates non-uniqueness, however, is not essential.
%Also, Tucker factorization generalizes on the CANDECOMP/PARAFAC (CP) decomposition, i.e., another popular tensor factorization method that decomposes the input tensor to a super-diagonal core tensor and factor matrices.
%In any tensor factorizations,
Each factor matrix $\mathbf{A}^{(n)}$ represents the latent features of the object related to the $n$th mode of $\T{X}$, and each element of a core tensor $\T{G}$ indicates the weights of the relations composed of columns of factor matrices.

However, in real-world, data are often incomplete with some missing entries.
To accommodate for the missing data, a partially observable Tucker factorization is needed.
%\begin{definition}[Partially observable Tucker factorization] \label{def:partial}
    Given a tensor $\T{X}$ $(\in \mathbb{R}^{I_{1} \times ... \times I_{N}})$
    with observable entries $\Omega$,
    the goal of \textit{partially observable Tucker factorization} of $\T{X}$ is to find factor matrices $\mathbf{A}^{(n)}$ $(\in \mathbb{R}^{I_{n} \times J_{n}}, n=1,\cdots,N)$ and a core tensor $\tensor{G}~(\in \mathbb{R}^{J_1 \times...\times J_N})$, which minimize the following loss: % function:
    \vspace{-1mm}
    \begin{equation}\label{eq:TF_PARTIAL_LF}
       L_F(\T{G},\mat{A}^{(1)},...,\mat{A}^{(N)}) = \sum_{\forall\alpha \in\Omega}{\left(\T{X}_{\alpha}-\sum_{\forall\beta \in\T{G}}\T{G}_{\beta} \prod_{n=1}^{N}a^{(n)}_{i_{n}j_n}\right)^{2}} +\lambda(\|\T{G}\|^2_{F} + \sum_{n=1}^{N}{{\| \mat{A}^{(n)} \|}^2_{F}}),
    \end{equation}
where $\alpha$ is an observable entry $(i_1,...,i_N)$ of input tensor $\T{X}$,
$\beta$ is an element $(j_1,...,j_N)$ of core tensor  $\T{G}$, and
$\lambda > 0$ is a regularization parameter.
%\qed
%\end{definition}
Note that the reconstruction error in Eq.~\eqref{eq:TF_PARTIAL_LF} depends only on the observable entries of $\T{X}$, and $L_{F}$ regularization is used in Eq.~\eqref{eq:TF_PARTIAL_LF} to prevent overfitting.
%For vectors, $L_{F}=L_{2}$, and $L_{2}$ regularization has been widely used in optimizations to prevent overfitting.

Tucker factorization often results in dense core and factor matrices. One of the approaches for sparsifying results is by including a sparsity constraint in the form of $L_1$ norm, a.k.a., Lasso, into the objective function.
%\subsection{Tucker Sparse Factorization Methods}\label{sec:prelim:tf}
%\begin{definition}[Partially observable Tucker factorization via sparse regularizer] \label{def:partial_l1}
Given a tensor $\T{X}$ with observable entries $\Omega$,
    the goal of \textit{partially observable Tucker factorization via sparse regularizer} of $\T{X}$ is to find factor matrices and a core tensor that minimize the following loss:
    %    $\alpha$ and $\beta$ are an entry of $\T{X}$ and $\T{G}$, respectively
    \begin{equation}    \label{eq:TF_PARTIAL_L1}
       L_{1}(\T{G},\mat{A}^{(1)},...,\mat{A}^{(N)}) =  \sum_{\forall\alpha \in\Omega}{\left(\T{X}_{\alpha}-\sum_{\forall\beta\in\T{G}}\T{G}_{\beta}\prod_{n=1}^{N}a^{(n)}_{i_{n}j_n}\right)^{2}} + \lambda(\|\T{G}\|_{1} + \sum_{n=1}^{N}{{\| \mat{A}^{(n)} \|}_{1}})
    \end{equation}
%    \qed
%\end{definition}
\noindent which changed the $L_2$ regularization term of Eq.~\eqref{eq:TF_PARTIAL_LF} to $L_1$.
Again the reconstruction error in Eq.~\eqref{eq:TF_PARTIAL_L1} depends only on the observable entries of $\T{X}$, and $L_{1}$ regularization is used to enforce sparsity.
Another approach for sparsifying results is by pruning. That is, a \textit{partially observable Tucker factorization via minimal element value pruning} of $\tensor{X}$ is to optimize on either Eq.~\eqref{eq:TF_PARTIAL_LF} or Eq.~\eqref{eq:TF_PARTIAL_L1}, and sets the smallest $s$ ratio of the elements to zero in the core and factor matrices.

Evaluation of tensor decomposition and the prediction of the missing entry values (a.k.a., tensor completion) involves reconstruction.
%\begin{definition}[Reconstruction] %, entry-wise reconstruction, reconstruction error, partial Reconstruction]
%\label{def:recon}
Given core tensor $\T{G}$ and factor matrices $\mathbf{A}^{(n)}$, the \textit{reconstruction} of the original tensor $\T{X}$ is defined as $ \T{X} \approx \T{G}\times_{1}\mathbf{A}^{(1)} \cdots \times_{N}\mathbf{A}^{(N)}$.

\subsection{Tucker ALS Algorithm}\label{sec:prelim:ALS}
%\begin{definition}[Alternating least squares]    \label{def:ALS}
A widely used technique for minimizing the loss functions Eq.~\eqref{eq:TF_PARTIAL_LF} and Eq.~\eqref{eq:TF_PARTIAL_L1} in a standard tensor factorization is alternating least squares (ALS)~\cite{kolda2009tensor}, which updates a factor matrix or a core tensor while keeping all others fixed.

    \begin{algorithm}[h!]
        %\small
        \caption{Tucker-ALS for Fully Observable Tensors (HOOI)} \label{alg:ALS_FULL}
        \SetKwInOut{Input}{Input}
        \SetKwInOut{Output}{Output}
        \Input{
            Tensor $\T{X} \in \mathbb{R}^{I_1 \times I_2 \times \cdots \times I_N}$, and
            core tensor dimensionality $J_1,...,J_N$. \\
        }
        \Output{
            Factor matrices $\mathbf{A}^{(n)} \in \mathbb{R}^{I_n \times J_n}$  $(n=1, ... , N)$, and
            core tensor $\T{G} \in \mathbb{R}^{J_1 \times J_2 \times \cdots \times J_N}$. \\
        }
        initialize all factor matrices $\mathbf{A}^{(n)}$\\
        \Repeat{reconstruction error converges or exceeds maximum iteration} {
            \For{$n=1...N$}{
                $\T{Y} \leftarrow \T{X}\times_{1}\mathbf{A}^{(1)\mat{T}} \cdots \times_{n-1}\mathbf{A}^{(n-1)\mat{T}}\times_{n+1}\mathbf{A}^{(n+1)\mat{T}}\cdots\times_{N}\mathbf{A}^{(N)\mat{T}} $ \\
                $\mathbf{A}^{(n)} \leftarrow J_{n}$ leading left singular vectors of $\T{Y}_{(n)}$
            }
        }
        $\T{G} \leftarrow \T{X}\times_{1}\mathbf{A}^{(1)\mathsf{T}} \cdots \times_{N}\mathbf{A}^{(N)\mathsf{T}}$
    \end{algorithm}

Algorithm~\ref{alg:ALS_FULL} describes a vanilla Tucker factorization algorithm based on ALS, which is called the \textit{higher-order orthogonal iteration} (HOOI) (see~\cite{kolda2009tensor} for details) that works on fully observable tensor.
Notice that Algorithm~\ref{alg:ALS_FULL} assumes missing entries of $\T{X}$ as zeros during the update process (lines 4-5).
However, setting missing values to zero enforces Tucker ALS to factorize the original tensor such that missing values becomes zero when reconstructed. Note that the missing values are often nonzero values that are unknown. Thus, setting missing values to zero inserts false information into the factorization which results in higher reconstruction error as well as higher generalization error.
%2 sparsity issue
Moreover, Algorithm~\ref{alg:ALS_FULL} computes SVD (singular vector decomposition) given $\T{Y}_{(n)}$, which often results in dense matrices, thus tensor-ALS results in overall dense core tensor and factor matrices.
%3 scalability issue - memory
Also, Algorithm~\ref{alg:ALS_FULL} requires storing a full-dense matrix $\T{Y}_{(n)}$, and the amount of memory needed for storing $\T{Y}_{(n)}$ is $O(I_n\prod_{m \neq n}{J_m})$.
The required memory grows rapidly when the order, the mode dimensionality, or the rank of a tensor increase, and ultimately causes \textit{intermediate data explosion}~\cite{kang2012gigatensor}.

In summary, the vanilla Tucker-ALS algorithm results in high generalization error in the presence of missing data, results in dense and thus hard-to-interpret core tensor and factor matrices, and cannot be applied to large data.
Therefore, Algorithm~\ref{alg:ALS_FULL} needs to be revised to focus only on observed entries, make sparse outputs, and be scaled for large-scale tensors at the same time.

\section{Proposed Method}
\label{sec:proposed_method}
In this section, we propose \vest (Very Sparse Tucker factorization), a method for partially observed large scale tensor, that results in very sparse core tensor and factor matrices.
Sparse results of tensor factorization increase interpretability and provide a scheme for better compression.
To maximize sparsity without losing accuracy, \vest iteratively updates core tensor and factor matrices, and prunes unimportant elements from the core tensor and factor matrices. However, there are several challenges in designing an efficient update and pruning rules.

\begin{itemize}[noitemsep, topsep=2pt]
\item \textbf{Evaluating importance of elements.} Vital elements of core tensor and factor matrices should not be pruned. How can we evaluate their importance?
%\item \red{\textbf{Keeping pruned elements to zero during update.} Pruned elements of core tensor and factor matrices should not be resurrected. How can we update efficiently without modifying pruned elements?}
%
\item \textbf{Automatically determining the sparsity.} There is a trade-off relationship between sparsity and accuracy. How can we automatically determine an appropriate sparsity which gives a good balance with regards to accuracy?
\item \textbf{Updating factors while guaranteeing non-decreasing sparsity.} The update process of factors and the core in the regular Tucker-ALS does not guarantee
that the sparsity improves over the update process. How can we guarantee that update rules improve the sparsity?
%
%\item \blue{\textbf{Non-decreasing sparsity.} In a regular Tucker-ALS, zero value elements can be updated to non-zeros. How can we ensure that update rules do not decrease the sparsity?}
\end{itemize}

\noindent We have the following main ideas to address the above challenges which we describe in detail in later subsections.
\begin{itemize}[noitemsep, topsep=2pt]
%\item \textbf{Suggest objective function} over that is dependent only on observable elements in the input tensor $\T{X}$.
\item \textbf{Design responsibility indicator} to evaluate contribution of each element on the accuracy
%Significant elements are conserved when truncating elements with low \textit{responsibility}.
(Section \ref{subsec:resp}).
\item {\textbf{Design auto-search algorithm \vesta} to find a good sparsity that resides near the maximum sparsity just before the reconstruction error shoots up} (Section~\ref{subsec:pruning}).
\item \textbf{Design element-wise update rules} to independently update each element of factor matrices and the core tensor. Element-wise update rules guarantee that the sparsity non-decreases by keeping pruned elements to zeros   (Sections \ref{subsec:update} and \ref{method:subsec:parallel}).
%We provide two kinds of element-wise update rules, derived from $L_F$ and $L_1$ regularization.
%
\end{itemize}

\subsection{Overview}

\begin{algorithm} [h!]
        %\small
        \caption{\vest: Very Sparse Tucker Factorization} \label{alg:overview}
        \SetKwInOut{Input}{Input}
        \SetKwInOut{Output}{Output}
        \Input{
            Tensor $\T{X} \in \mathbb{R}^{I_1 \times I_2 \times \cdots \times I_N}$,
            core tensor dimensionality $J_1,...,J_N$, and
            target sparsity $s$ (if manual-mode).
            %initial pruning rate $init\_pr$,
            %and maximum pruning rate $max\_pr$
        }
        \Output{
            Sparse factor matrices $A^{(n)} \in \mathbb{R}^{I_n \times J_n} (n=1, \cdots, N)$ and
            sparse core tensor $\T{G}$ $\in \mathbb{R}^{J_1 \times J_2 \times \cdots \times J_N}$. \\
        }
        \vspace{1.5mm}
        randomly initialize $\mathbf{A}^{(n)} (n=1,\cdots,N)$ and $\T{G}$; set $pr$ = {\scriptsize INIT$\_$PR}, $iterN$ = $0$.\\
        \Repeat{RE converges or iterN++ exceeds maximum iteration} {
            update unpruned elements of $A^{(n)} (n=1,\cdots,N)$ \Comment*[f]{ Algorithm~\ref{alg:FMUpdate}}\\
            update unpruned elements of $\T{G}$  \Comment*[f]{ Eq.\eqref{eq:coreupdate_lf} or \eqref{eq:coreupdate_l1}} \\
            compute RE using observable entries $\Omega$ \Comment*[f]{Eq.\eqref{eq:re}} \\
            \If{should\_prune()}{
                prune $pr$=min({\scriptsize INIT$\_$PR}*$iterN$, {\scriptsize MAX$\_$PR}) ratio of elements $e$ in $A^{(n)}$ and $\T{G}$ based on $Resp(e)$ of the elements \Comment*[f]{Algorithm~\ref{alg:Pruning}}\\
            }
        }
    \For{n = 1...$N$}{
        $ \mathbf{U}^{(n)}\mathbf{B}^{(n)} \leftarrow \mathbf{A}^{(n)}$, and set $\mathbf{A}^{(n)} \leftarrow \mathbf{U}^{(n)}$\\
        $\tensor{G} \leftarrow \tensor{G}\times_{n}\mathbf{B}^{(n)}$
    }
\end{algorithm}

\vest is a scalable Tucker factorization method that results in very sparse core tensor and factor matrices for partially observable data (see Algorithm~\ref{alg:overview}).
First, \vest initializes all elements of the core tensor and factor matrices with random real values between 0 and 1 (line 1).
Next, \vest iteratively updates the core tensor and factor matrices while pruning their elements (lines 3-7).
In lines 3-4, \vest updates unpruned elements of the core tensor and factor matrices by element-wise update rules (Section~\ref{subsec:update}),
guaranteeing that the sparsity non-decreases.
Then \vest prunes unimportant elements in the core tensor and factor matrices (lines 6-7).
Importance of each element $e$ is evaluated by responsibility $Resp(e)$ which indicates how largely the element contributes to the accuracy (Section~\ref{subsec:resp}).
$should\_prune()$ function determines when to stop pruning: if desired sparsity $s$ is achieved (in the manual version \vestm) or the reconstruction error shows a rapid increase (in the automatic version \vesta).
Motivated from simulated annealing, we gradually increase the pruning ratio as iterations proceed (line 7); this enables to explore larger search space in the beginning, while reducing the extent of the search to reduce to a minimum in the later iterations.
%
%In the case of \vesta, it revives the latest pruned elements to recover Accuracy (in step 3) explained in \ref{subsec:auto}.
%
The iterations proceed until the reconstruction error converges or the maximum iteration is reached.
Finally, \vest standardizes all columns of factor matrices such that their norm is equal to one, and updates core tensor accordingly (lines 9-11). Specifically, $A^{(n)}$ is decomposed to $\mathbf{U^{(n)}} \in \mathbb{R}^{I_n\times J_n}$ and $\mathbf{B^{(n)}}$ where columns of $\mathbf{U^{(n)}}$ are unit vectors and $\mathbf{B^{(n)}}$ is a diagonal matrix whose $(i, i)$th element is the norm of $A^{(n)}$'s $i$th column. The core tensors are updated to maintain the same reconstruction error \cite{kolda2006multilinear}.

\subsection{Evaluating Importance of Elements by \textit{Responsibility}} \label{subsec:resp}
\label{sec:method:object}

\vest calculates the \textit{responsibility} of each element which represents its contribution to the overall reconstruction accuracy over the observable elements of the input tensor to determine and prune unimportant elements.
The intuition is that reconstruction error increases significantly when a vital element of the core tensor or factor matrices is set to zero, i.e., pruned.
On the other hand, if the reconstruction error after pruning is similar or even smaller to that before pruning, the pruned element is insignificant.
Formally, the responsibility is defined as follows.
\begin{definition}[Responsibility] \label{def:res_base}
    Responsibility of an element $e$ in a factor matrix ($e=a^{(n)}_{ij}$) or core tensor ($e=\T{G}_\gamma$) is given by
   % ${Resp(e)} = \frac{RE(e)- RE}{RE}$,
    \begin{equation}\label{eq:resp}
    {Resp(e)} = \frac{RE(e)- RE}{RE},
    \end{equation}
    where
    \begin{equation}\label{eq:re}
    RE = \sqrt{\sum\limits_{\forall\alpha=(i_1,\cdots, i_N)\in\Omega} \left(\tensor{X}_{\alpha} - {\sum\limits_{\forall{\beta=(j_1,\cdots, j_N)\in\T{G}}}\tensor{G}_{\beta}\prod\limits_{n=1}^{N}a^{(n)}_{i_nj_n}}\right)^2}/{||\tensor{X}||_F}
    \end{equation}
    is the normalized reconstruction error over the observable entries $\Omega$ of the original tensor $\T{X}$,
    and $RE(e)$ is residual reconstruction error defined when the element $e$ is set to zero (Eq.~\eqref{eq:resp_core} and \eqref{eq:resp_fac}).
\qed
\end{definition}

\begin{definition}[Residual reconstruction error] \label{def:res_RRE}
\noindent
The residual reconstruction error $RE(\T{G}_{\gamma})$ for $(j_1,...,j_N)$th element $\gamma$ in core tensor $\T{G}$ is as follows:
    % reconstruction error for core matrix
%    \red{update needed}
    \begin{equation}\label{eq:resp_core}
    (RE(\T{G}_{\gamma}))^2 = \frac{\sum\limits_{\forall\alpha=(i_1,\cdots, i_N) \in \Omega} \big( \bm{\mathscr{X}}_\alpha-\sum\limits_{\forall\beta=(j_1,\cdots, j_N) \neq\gamma\in\T{G}}\bm{\mathscr{G}}_\beta \prod\limits_{n=1}^{N}a^{(n)}_{i_{n}j_{n}}\big)^{2}}{||\tensor{X}||_F^2}
    \end{equation}

\noindent The residual reconstruction error $RE(a^{(n)}_{i,j})$ for an $(i, j)$th element $a^{(n)}_{i,j}$ in a factor matrix $\mathbf{A}^{(n)}$ is as follows:
\begin{equation} \label{eq:resp_fac}
\begin{split}
(RE(a^{(n)}_{ij}))^2 =  & %\frac{\sum\limits_{\forall\alpha\in\Omega_{i_n}^{(n)}}\left(\T{X}_{\alpha} -B_{j_n\neq j}(\alpha)\right)^2 + \sum\limits_{\forall\alpha\notin\Omega_{i_n}^{(n)}}\left(\T{X}_{\alpha}-B(\alpha)\right)^2}{||\tensor{X}||_F^2}\\
%= &
{RE}^2 + \frac{\sum\limits_{\forall\alpha\in\Omega_{i_n}^{(n)}} \big(2\cdot(\T{X}_{\alpha} -B(\alpha)) + B_{j_n=j}(\alpha)\big)\cdot(B_{j_n=j}(\alpha))}{||\tensor{X}||_F^2},
\end{split}
\end{equation}
where $B(\alpha)$ is the \textit{entry-wise reconstruction} defined as
\begin{equation} \label{eq:recon2}
    \tensor{X}_{\alpha=(i_1,\cdots, i_N)} \approx B(\alpha) =
        \sum\limits_{\forall{\beta=(j_1,\cdots, j_N)}\in\T{G}}\tensor{G}_{\beta}\prod\limits_{n=1}^{N}a^{(n)}_{i_nj_n},
\end{equation}
and $B_{j_n=j}(\alpha)$ and $B_{j_n\neq j}(\alpha)$ are the \textit{partial reconstruction functions} defined as
\begin{equation} \label{eq:precon}
    B_{j_n=j}(\alpha) =   \sum\limits_{\forall\beta_{j_n=j}}\T{G}_{\beta_{j_n= j}} \prod\limits_{n=1}^{N}a_{i_nj_n}^{(n)},\;\;
    B_{j_n\neq j}(\alpha) =  \sum\limits_{\forall\beta_{j_n\neq j}}\T{G}_{\beta_{j_n\neq j}} \prod\limits_{n=1}^{N}a_{i_nj_n}^{(n)}.
\end{equation}
\noindent
\qed
\end{definition}
% Notice that $B(\alpha) =  B_{j_n=j}(\alpha) + B_{j_n\neq j}(\alpha)$.
% Eq.~\eqref{eq:recon2} is used to predict values of missing entries.

The proof of correctness for the derivation of the Eq.\eqref{eq:resp_fac} is provided in the supplementary material~\cite{supp-material}.
Note that both definitions are derived from the element-wise reformulation of the reconstruction error.

\subsection{Pruning}
\label{subsec:pruning}
After calculation of responsibility values, \vest prunes core tensor and factor matrices.
Pruning is performed iteratively, each time after the core tensor and factor matrices are updated.
The process of pruning an element consists of setting the value of the element to zero and marking the element as pruned in a marking table. The marked elements are excluded from the update step.
To prune elements with low responsibility values, \vest sorts elements of the core tensor and each factor matrix, respectively, by the responsibility $Resp(e)$ in ascending order.
Then, \vest prunes smallest $pr|\T{G}|$ elements from core tensor and smallest $pr|A^{(n)}|$ from each factor matrix, where $pr$ is the pruning rate of the current iteration.
\vest starts with a small pruning rate $pr$ (INIT$\_$PR) and slowly increases $pr$ until maximum pruning rate (MAX$\_$PR) is reached.
The default values of INIT$\_$PR and MAX$\_$PR are set to $0.01$ and $0.1$, respectively.

To determine when to stop pruning, \vest provides two different algorithms, the manual version \vestm and the automatic version \vesta; the subscript * denotes whether $L1$ or $L_F$ regularization is used.
\begin{itemize}[noitemsep,topsep=0pt]
    \item \vestm takes a target sparsity $s$ as an input from the user and stops pruning when the total sparsity reaches $s$. That is, \vestm enables users to decide on the lower bound of the final sparsity.
    \item \vesta determines the final sparsity automatically.
    \vesta does this by tracking changes in the reconstruction error and determines to stop pruning when an elbow point of the reconstruction error curve is reached. 
    The elbow point is estimated as the point when the second derivative of the RE curve, estimated as $(RE_{t} + RE_{{t-2}} - 2*RE_{{t-1}})/pr_{t}$ where $RE_{t}$ and $pr_t$ are RE and pruning rate at $t^{th}$ iteration, respectively, exceeds a small threshold (0.05 used).
\end{itemize}

%For \vesta, the elbow point is estimated as the point when the second derivative of the RE curve, estimated as $(RE_{t} + RE_{{t-2}} - 2*RE_{{t-1}})/pr_{t}$ where $RE_{t}$ and $pr_t$ are RE and pruning rate at $t^{th}$ iteration, exceeds a small threshold (0.05 used).

\subsection{Element-Wise Update Rules}\label{subsec:update}
\vest updates elements of the core tensor and factor matrices based on a coordinate decent  approach in parallel.
It enables \vest to update the core tensor and factor matrices without changing the value of the pruned elements.
\vest checks the marking table which indicates whether elements have been pruned, and updates only the un-pruned elements.
The update of an element is performed with observable tensor entries and fixed values of other elements in the factor matrices and the core tensor.
The update rules for the core tensor and factor matrices are derived by setting the partial derivative of the loss function to zero and solving for each element.
In previous works~\cite{Oh2018,Lee2018a}, this approach has been shown and proven to converge faster with higher accuracy than existing approaches. %,Oh2019}.
Advantages of our update rules are that
1{)} accuracy is high and convergence is faster~\cite{Oh2018},
2{)} parallelization and selective updates are possible because all the elements are independently updated, and
3{)} the size of intermediate data is small, making the algorithm scalable.

\subsubsection{Element-wise update rules with $L_F$ regularization}
%We derive a partial differential equation for updating factor matrices.
The update rule for an element $a_{i_nj_n}^{(n)}$ of factor matrix $A^{(n)}$ is derived by setting the partial derivative of loss function (Eq.~\eqref{eq:TF_PARTIAL_LF}) with regard to $a_{i_nj_n}^{(n)}$ to zero.
%Lemma~\ref{lemma:FMupdate_lff} shows how to update $a_{i_nj_n}^{(n)}$ of $\mathbf{A}^{(n)}$ using the partial derivative.

\begin{lemma}[Update rule for factor matrices with $L_F$ regularization]\label{lemma:FMupdate_lff}
    \begin{equation}\label{eq:factor_lf}
        a_{i_nj_n}^{(n)} \longleftarrow  \argmin{a_{i_nj_n}^{(n)}}{\mathcal{L}_F(\T{G},\mat{A}^{(1)},...,\mat{A}^{(N)})} =
        \frac{ \big( \sum\limits_{\forall\alpha\in\Omega_{i_n}^{(n)}}\bm{\mathscr{X}}_{\alpha}\delta_{\alpha}^{(n)}(j_n) \big) - \big( \sum\limits_{\forall t\neq j_n}\bm{v}_{i_nj_n}^{(n)}(t)\cdot a_{i_nt}^{(n)} \big) }
{\bm{v}_{i_nj_n}^{(n)}(j_n) + \lambda},
    \end{equation}
    where $\bm{v}_{i_nj_n}^{(n)}$ is a length $J_n$ vector whose $j$th element is
   \begin{equation}\label{eq:FMU_v}
\bm{v}_{i_nj_n}^{(n)}(j) = \sum\limits_{\forall\alpha\in\Omega_{i_n}^{(n)}}\delta_{\alpha}^{(n)}(j)\delta_{\alpha}^{(n)}(j_n),
   \end{equation}
    $\delta_\alpha^{(n)}$ is a length $J_n$ vector whose $j$th element is
\begin{equation}\label{eq:FMU_delta}
 \delta_\alpha^{(n)}(j) = \sum\limits_{\forall\beta_{j_n=j}\in\T{G}} \T{G}_{\beta_{j_n=j}}\prod\limits_{k\neq n} a_{i_kj_k}^{(k)},
    \end{equation}
    $\Omega_{i_n}^{(n)}$ is the subset of $\Omega$ whose index of $n$th mode is $i_n$, and $\lambda$ is a regularization parameter. \qed
    \end{lemma}

The derivation of the core tensor update rule is similar to that of the factor matrix.
The update rule for the $\beta^{th}$ element $\T{G}_{\beta}$ of the core tensor $\T{G}$ is given as follows.
\begin{lemma}[Update rule for core tensor with $L_F$ regularization]\label{lemma:FMupdate_lfc}
    %\small
    \begin{equation} \label{eq:coreupdate_lf}
    \bm{\mathscr{G}}_{\beta} \longleftarrow \frac{\sum\limits_{\forall\alpha\in\Omega}(\bm{\mathscr{X}}_\alpha-\sum\limits_{\forall\gamma\neq\beta}\bm{\mathscr{G}}_\gamma \prod\limits_{n=1}^{N}a^{(n)}_{i_{n}j_{n}})\cdot \prod\limits_{n=1}^{N}a^{(n)}_{i_{n}j_{n}}}{\lambda + \sum\limits_{\forall\alpha\in\Omega}(\prod\limits_{n=1}^{N}a^{(n)}_{i_{n}j_{n}})^2}
    \end{equation}
    \qed
\end{lemma}

\subsubsection{Element-wise update rules with $L_1$ regularization}
    For an element $a_{i_nj_n}^{(n)}$ of factor matrix $A^{(n)}$, the element-wise update rule with $L_1$ regularization is provided in the following Lemmas.
\begin{lemma}[Update rule for factor matrix with $L_1$ regularization]\label{lemma:FMupdate_l1f}
    %\small
    \begin{equation}\label{eq:F1}
    \argmin{a_{i_nj_n}^{(n)}}{L_{1}(\tensor{G},\mat{A}^{(1)},...,\mat{A}^{(N)})} =
        \begin{cases}
        (\lambda - g_{fm})/d_{fm} & \text{if $g_{fm} > \lambda$} \\
        -(\lambda + g_{fm})/d_{fm} & \text{if $g_{fm} < -\lambda$} \\
        0 & \text{otherwise}
        \end{cases}
    \end{equation}
    where
%\begin{equation}\label{eq:FM-g}
$g_{fm} = 2 \big( \sum\limits_{\forall\alpha\in\Omega_{i_n}^{(n)}}\bm{\mathscr{X}}_{\alpha}\delta_{\alpha}^{(n)}(j_n) \big) - \big( \sum\limits_{\forall t\neq j_n}\bm{v}_{i_nj_n}^{(n)}(t)\cdot a_{i_nt}^{(n)} \big)$,  %\end{equation}
%\begin{equation}\label{eq:FM-d}
$d_{fm} = 2\bm{v}_{i_nj_n}^{(n)}(j_n)$, and
%\end{equation}
$\bm{v}_{i_nj_n}^{(n)}$,  $\delta_\alpha^{(n)}$,
$\Omega_{i_n}^{(n)}$, and $\lambda$ follow the same specification provided in Lemma~\ref{lemma:FMupdate_lff}. \qed
\end{lemma}

For an element $\tensor{G}_\beta$ of core tensor, the element-wise update rule with $L_1$ regularization is as follows:
\begin{lemma}[Update rule for core tensor with $L_1$ regularization]\label{lemma:FMupdate_l1c}
    %\small
    \begin{equation}\label{eq:coreupdate_l1}
        \argmin{\tensor{G}_\beta}{L_{1}(\tensor{G},\mat{A}^{(1)},...,\mat{A}^{(N)})} =
        \begin{cases}
            (\lambda - g_c)/d_c & \text{if $g > \lambda$} \\
            -(\lambda + g_c)/d_c & \text{if $g < -\lambda$} \\
            0 & \text{otherwise}
        \end{cases}
    \end{equation}
  where
%  \begin{equation}\label{eq:core-g}
  $g_c = -2\sum\limits_{\forall\alpha\in\Omega}\big(\tensor{X}_\alpha-\sum\limits_{\forall\gamma\neq\beta}\tensor{G}_\gamma \prod\limits_{n=1}^N a_{i_nj_n}^{(n)}\big)\cdot\prod\limits_{n=1}^Na_{i_nj_n}^{(n)}$, %\end{equation}
  and
%  \begin{equation}\label{eq:core-d}
  $d_c = 2\sum\limits_{\forall\alpha\in\Omega}\big(\prod\limits_{n=1}^N a_{i_nj_n}^{(n)}\big)^2$.
%  \end{equation}
  \qed
\end{lemma}

The proofs of Lemmas~\ref{lemma:FMupdate_lff} to~\ref{lemma:FMupdate_l1c} are provided in the supplementary material~\cite{supp-material}.

\subsection{Parallel Update Algorithms}
\label{method:subsec:parallel}
Responsibility calculation and factor matrices updates are performed in parallel.
Algorithm~\ref{alg:Pruning} describes the pruning process where responsibility values of the core tensor and factor matrices are calculated in parallel for each observable entries of the input tensor.
Note that the use of $B(\alpha)$ in line 5 enabled fast computing of $Resp(\T{G}_\beta)$ in line 6; for a given $\beta$ in line 4, computing line 5 requires $O(|\Omega|)$ rather than $O(|\Omega| |\tensor{G}|)$ since there is no need to compute $\sum_{\forall\beta\neq\gamma\in\T{G}}{\T{G}}_{\gamma} \prod_{n=1}^{N}a^{(n)}_{i_{n}j_{n}}$ in Eq.~\eqref{eq:resp_core} from scratch.

\begin{algorithm} [h!]
    %\small
    \caption{Parallel Pruning} \label{alg:Pruning}
    \SetKwInOut{Input}{Input}
    \SetKwInOut{Output}{Output}
    \Input{
        Tensor $\T{X} \in \mathbb{R}^{I_1 \times \cdots \times I_N}$,
        factor matrices $\mathbf{A}^{(n)} \in \mathbb{R}^{I_n \times J_n} (n=1,\cdots,N)$,
        core tensor $\T{G} \in J_1 \times ... \times J_N$, and
        pruning ratio $pr$.
    }
    \Output{
        Pruned $\mathbf{A}^{(n)} (n=1,\cdots,N)$ and $\T{G}$
    }
    \vspace{1.5mm}
    \For(\Comment*[f]{\textbf{in parallel}}){$\alpha = \forall(i_1,...,i_N)\in\Omega$}
        {
            calculate $B(\alpha) = \sum\limits_{\forall{\beta=(j_1,\cdots, j_N)}\in\T{G}}\tensor{G}_{\beta}\prod\limits_{n=1}^{N}a^{(n)}_{i_nj_n}$ \\
            calculate $\T{X}_\alpha - B(\alpha)$  \Comment*[f]{Eq.\eqref{eq:re}}    \\
            % \For{$\beta = \forall(j_1,...,j_N)\in\T{G}$}
            % {
            %     Calculate $\T{X}_\alpha - B(\alpha)$ \\
            % }
        }
        \For(\Comment*[f]{\textbf{in parallel}}){$\beta = \forall(j_1,...,j_N)\in\T{G}$}
       {
            calculate $\sum\limits_{\forall\alpha\in\Omega}( \bm{\mathscr{X}}_\alpha - B(\alpha) + \bm{\mathscr{G}}_\beta \prod\limits_{n=1}^{N}a^{(n)}_{i_{n}j_{n}} )$ \\
            % \For{$\alpha = \forall(i_1,...,i_N)\in\Omega$}
            % {
            %     Calculate $\prod_{n=1}^{N}a^{(n)}_{i_{n}j_{n}}$ \\
            % }
            calculate $Resp(\T{G}_\beta)$ \Comment*[f]{Eq.\eqref{eq:resp}, \eqref{eq:resp_core}} \\
        }

    sort core tensor elements by $Resp(\T{G}_\beta)$ values in an ascending order \\
    \For(){$i_n = 1 ... I_n$}
    {
        \For(\Comment*[f]{\textbf{in parallel}}){$j_n = 1 ... J_n$}
        {
            \For{$\alpha =     \forall(i_1,...,i_N)\in\Omega_{i_n}^{(n)}$}
            {
                calculate $(2(\T{X}_\alpha - B(\alpha))+B_{j_n=j}(\alpha))\cdot B_{j_n=j}(\alpha)$\\
            }
            calculate $Resp(a_{i_nj_n}^{(n)})$ \Comment*[f]{Eq.\eqref{eq:resp}, \eqref{eq:resp_fac}}     \\
        }
        sort factor matrix elements by $Resp(a_{ij}^{(n)})$ values in an ascending order \\
    }
    prune smallest $pr|\T{G}|$ and $pr|A^{(n)}|$ elements of $\T{G}$ and $A^{(n)} (n=1,...,N)$, respectively. \\
\end{algorithm}

\begin{algorithm} [ht!]
    %\small
    \caption{Parallel Element-Wise Factor Matrix Update} \label{alg:FMUpdate}
    \SetKwInOut{Input}{Input}
    \SetKwInOut{Output}{Output}
    \Input{
        Tensor $\T{X} \in \mathbb{R}^{I_1 \times \cdots \times I_N}$,
        factor matrices $A^{(n)} \in \mathbb{R}^{I_n \times J_n} (n=1,\cdots,N)$, and
        core tensor $\T{G} \in J_1 \times ...J_N$.
    }
    \Output{
        Updated factor matrices $A^{(n)} \in \mathbb{R}^{I_n \times J_n} (n=1,\cdots,N)$
    }
    \vspace{1.5mm}
    \For(\Comment*[f]{$n$th factor matrix}){$n = 1 ... N$}
    {
        \For(\Comment*[f]{\textbf{in parallel}}){$i_n = 1 ... I_n$}
        {
            \For{$j_n = 1 ... J_n$}
            {
                \If{ $a_{i_nj_n}^{(n)}$ is pruned} {continue}
                \For{$\alpha = \forall(i_1,...,i_N)\in\Omega_{i_n}^{(n)}$}
                {
                    \For(\Comment*[f]{compute $\delta$}){$\beta = \forall(j_1,...,j_N)\in\T{G}$}
                    {
                        $\delta_\alpha^{(n)}(j_n) \longleftarrow \delta_\alpha^{(n)}(j_n) + \T{G}_\beta\prod\limits_{\forall k \neq n}a_{i_kj_k}^{(k)}$\\
                    }
%                    accumulate $\sum\limits_{\forall\alpha\in\Omega_{i_n}^{(n)}}\T{X}_\alpha\delta_\alpha^{(n)}(j_n)$, and update $\bm{v}_{i_nj_n}^{(n)}$\Comment*[f]{Eq.\eqref{eq:FMU_v},\eqref{eq:FMU_delta}}\\
                    accumulate $\T{X}_\alpha\delta_\alpha^{(n)}(j_n)$, and update $\bm{v}_{i_nj_n}^{(n)}$\Comment*[f]{Eq.\eqref{eq:FMU_v},\eqref{eq:FMU_delta}}\\
                }

                update $a_{i_nj_n}^{(n)}$ using Eq. \eqref{eq:factor_lf} for $L_F$ (use Eq. \eqref{eq:F1} for $L_1$)
            }
        }
    }
\end{algorithm}

The element-wise update of factor matrix $A^{(n)}$ is performed in parallel for each rows of factor matrices using either the $L_F$ or $L_1$ regularization (see Algorithm~\ref{alg:FMUpdate}).
Elements of the core tensor are dependent on each other and thus cannot be updated in parallel.
However, considering that typical size ${|\T{G}|}$ of the core tensor is small, the core tensor updates are a minor burden in the computational process.

\begin{lemma}[Complexity of \vest per iteration]\label{lemma:complexity}\\
%\small
%Overall time complexity of \vest per iteration is $~O(N^2J|\T{G}||\Omega|/T)~$
The time complexity per iteration of \vest is $O(N^2J|\T{G}||\Omega|/T + N IJ \log (IJ) + |\T{G}| \log |\T{G}|)$,
and  the memory complexity is $~O(TJ+J^N+NIJ)~$,
where T is the number of threads, N is the order, I is the dimensionality of input tensor, and J is the dimensionality of the core tensor assuming the dimensionalities are equal for all modes.   \qed
\end{lemma}
The full proof of Lemma~\ref{lemma:complexity} is in the supplementary material~\cite{supp-material}.

\section{Experiments}
\label{sec:experiment}
We conduct experiments to answer the following questions.
\begin{enumerate}[noitemsep, topsep=2pt]
    \item \textbf{Performance comparison (Section~\ref{sec:accuracy}).} How accurately and sparsely does \vest decompose a given tensor compared to other methods? %How accurate does \vest predict missing entry values?
    \item \textbf{Sparsity and accuracy (Section~\ref{sec:result_sparsity}).}
    Does \vest successfully prune redundant information in the decomposition without hurting the accuracy?
    Does \vesta find reasonable sparsity and accuracy trade-off point?
    %What is the effect of different $\lambda$ values on the sparsity of \vestl model?
    \item \textbf{Data scalability (Section~\ref{sec:scalablility}).} How scalable is \vest?
    %How large is a tensor that \vest decompose?
%    \item \textbf{Factorization stability.} How stable are the \vest results?
    \item \textbf{Interpretable discoveries (Section~\ref{sec:discovery}).} How interpretable are the \vest results % compared to dense results
    for discoveries on real-world tensors?
\end{enumerate}

\subsection{Experimental Settings}\label{sec:experimental-setting}
\textbf{Datasets.}
We used three real-world datasets and synthetic datasets as summarized in Table~\ref{table:datasets}.
The real-world datasets are
MovieLens\footnote{\url{https://grouplens.org/datasets/movielens/}},
Yelp\footnote{\url{http://www.yelp.com/dataset\_challenge/}}, and
AmazonFood\footnote{ \url{http://snap.stanford.edu/data/web-FineFoods.html}}.
MovieLens is a $4^{th}$ order tensor of movie ratings containing (user, movie, year, hour).
Yelp is a $3^{rd}$ order tensor of business services rating data containing (user, business, year-month).
AmazonFood is a $3^{rd}$ order tensor of food review scores from Amazon containing (product, user, year-month).
To compare with other methods, we used subsets Yelp-s and AmazonFood-s of 3$^{rd}$ order tensors which are made denser than their originals.
The density of Yelp-s and AmazonFood-s are 0.01 and 0.02, respectively.
We also generated synthetic random tensors of various sizes and orders to test data scalability.

\begin{table}
    \centering
    \caption{Summary of datasets and hyperparameters used.}
\resizebox{\linewidth}{!}{
    \begin{tabular}{lcll llr}
        \toprule
        Name        &   Order   &  Dimensionality   & Ranks & $|\Omega|$ & $|\Omega|_{test}$ \\ % & $\lambda'_{L_F}$($\lambda'_{L_1}$)\\
        \midrule
        MovieLens   &     4     & $138K \times 27K \times 21 \times 24$ &$6\times 6 \times 2 \times 2$  &   18M  &  2M \\% & 5 \\ %90.47 \\
        Yelp        &     3     &     $71K  \times 16K \times 108$   & $10  \times 10 \times 10$   &  301K  &  33K \\% & 5 \\ %1.74\\
        AmazonFood  & 3 & $74K \times 256K  \times 143$ & $9  \times 9 \times 14$ & 511K & 57K \\% & 5 \\ % 0.86 \\
        Yelp-s        &     3     &     $50  \times 50 \times 10$   & $5  \times 5 \times 5$   &  235  &  32\\% & 5 \\ % 1.740741\\
        AmazonFood-s  & 3 & $50 \times 50  \times 10$ & $5  \times 5 \times 5$ & 444 & 51\\% & 5 \\ %3.288889\\
   %     Syn10/50/80 & 3 & $50 \times 50 \times 10$ & $5 \times 5 \times 5$ &
   %     2K/11K/18K & 2H/1K/2K \\ %& 5 \\
        Synthetic   & $3-10$ &              $10^3-10^8$             & $3 \times \cdots  \times 3$ & $10^3-10^7$  &  - \\% \\ & 5\\
        \bottomrule
    \end{tabular}
    }
    \label{table:datasets}
\end{table}

\noindent \textbf{Environment.}
\vest was written in C++ with \omp \cite{Dagum:1998:OIA:615255.615542} and \arma \cite{sanderson2016armadillo} libraries for parallelization.
Methods L1 (Lasso) and Value Pruning were run on \vest framework with the difference just in the pruning approaches.
We used the codes provided by the authors for TTP~\cite{Sun2017} (R) and Sparse CP~\cite{allen2012sparse} (Matlab).
Tucker-ALS was performed via Tensor Toolbox for Matlab~\cite{TTB_Software}.
All experiments were done on a single machine equipped with an Intel Xeon E5-2630 v4 2.2GHz CPU (10 cores/20 threads) and 512GB memory.
All reported measures are averages of five runs, unless otherwise stated.

\noindent \textbf{Competitors.}
We compared \vest with the following methods.
\begin{enumerate}[noitemsep,topsep=2pt]
%    \item[$\bullet$] %CP-ALS~\cite{acar2010scalable}: CP decomposition for missing entry prediction.
    \item[$\bullet$] L1 (Lasso): A Tucker factorization method with lasso sparsity constraint implemented as \vestml with sparsity $s=0$.
    \item[$\bullet$] Value Pruning: A Tucker factorization method with value pruning at the last step implemented as \vestm with sparsity $s=0$ followed by value pruning with ratio 0.6 for $L_F$ and $L_1$ losses.
    \item[$\bullet$] TTP~\cite{Sun2017}: A tensor decomposition method that results in sparse components.
    \item[$\bullet$] Sparse CP~\cite{allen2012sparse}: CP decomposition method with lasso penalty.
    \item[$\bullet$] Tucker-ALS~\cite{kolda2009tensor}: Conventional Tucker factorization method (HOOI).
%    \item[$\bullet$] Tucker-ALS with pruning: Tucker-ALS combined with pruning of small cell values according to Definition~\ref{def:partial_minFV}.
    \item[$\bullet$] PTucker-Approx~\cite{Oh2018} : Tucker decomposition method for partially observable tensor with iterative element value pruning.

\end{enumerate}

\subsection{Performance Comparison}\label{sec:accuracy}
We compared the accuracy of \vestaf and \vestal with those of the competitors on datasets Yelp-s and AmazonFood-s (Table~\ref{table:datasets}).
The comparison was performed on smaller and denser datasets of order three and not on original real-world datasets due to limitations of the TTR and Sparse CP.
%Details of hyper-parameter selections are provided in the supplementary material~\cite{supp-material}.

We measured and compared normalized reconstruction errors (RE) over observable entries in input tensors.
As shown in Fig.~\ref{fig:accuracy}\textbf{(a)}, \vestaf and \vestal decomposed a given tensor with at least 2.8 times lower RE compared to other methods at a similar sparsity. % for partially observable tensor.
Fig.~\ref{fig:accuracy}\textbf{(a)} also shows that at a similar RE value, \vestf and \vestl output at least 2.2 times more sparse factor matrices and core tensor compared to other methods, where the sparsity is measured as the ratio of number of nonzero values in $\T{G}$ and $\mathbf{A}^{(n)}$ over $|\T{G}| + |\sum_{n=1}^{N}{\mathbf{A}^{(n)}}|$.

To answer how well \vest predicts missing entries compared to other methods, we measured REs of the reconstructed missing values (Test RE).
After learning the factor matrices and the core tensor using 90\% of the observed entries, we calculated the Test REs on the remaining $10\%$ of the observable entries.
Fig~\ref{fig:accuracy}\textbf{(b)} shows that \vest predicts missing entries at least 1.8 times more accurately compared to others, in addition to providing at least 1.7 times sparser results.

%In both reconstruction and missing entry tests, the low accuracies of Tucker-ALS and CP-ALS are due to the assumption that values of unknown or missing entries are set to zero, which enforces the factorization to optimize toward the zero values.
%The difference in the reconstruction errors becomes larger as the input data becomes more sparse.

% Note that \red{when entries are truncated according to the magnitude of entry values, the accuracy degrades significantly.
% We also observe that CP-ALS is not very sparse when sparsity of factor matrix is considered even though it learns a very sparse core tensor.}

\subsection{Sparsity and Accuracy}
\label{sec:result_sparsity}

%\red{PUT fig\ref{fig:sparsity_recon} in the Supplemental Materials .}
\begin{figure}[h!]
    \centering
    \includegraphics[width=.9\textwidth]{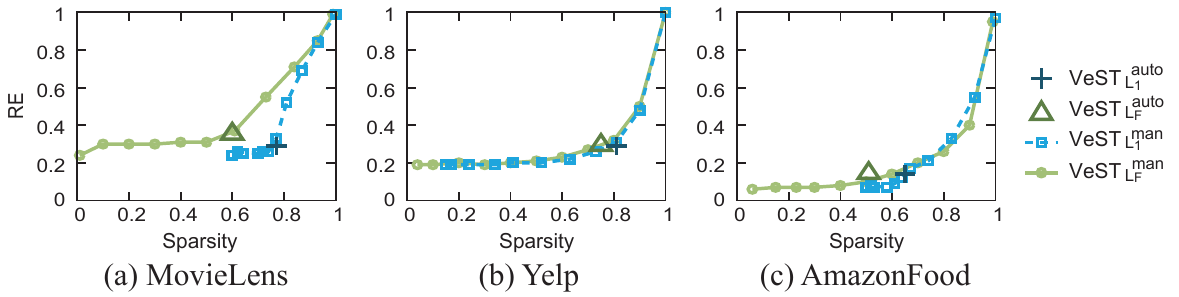}
    \caption{Sparsity against reconstruction error (RE) of \vest with varying target sparsity $s$.
    Note that
    1) \vest successfully removes redundant information in the decomposition without hurting the accuracy,
    and
    2) \vesta automatically finds the sparsity at the elbow point of the RE curve.
%    Measurements are averages of five runs.
}
    \label{fig:sparsity_RE}
\end{figure}

We tested the sparsity and accuracy of outputs of \vest on three full size real-world datasets: MovieLens, Yelp, and AmazonFood, in Figure~\ref{fig:sparsity_RE}. % as shown in Table~\ref{table:datasets}.
%\textbf{sparsity v.s. Accuracy.}
%
First, we investigated how the sparsity affects the normalized reconstruction error (RE).
Note that the REs are not affected much until the sparsities are above 0.6; this shows that there is redundant information in the decomposition results, and \vest successfully finds and removes such redundancies to get compact outputs.
% We investigated how sparse \vest makes the core tensor and the factor matrices by changing the target sparsity $s$ of \vestmf and \vestml, and evaluating their sparsity and accuracy measures.
% Fig.~\ref{fig:sparsity_RE} shows the overall sparsity against the normalized reconstruction error (RE) when the target sparsity $s$ is varied from $0$ to $0.99$.
%
%
Second, we investigated how \vesta automatically finds a desired sparsity which gives a good tradeoff with regard to accuracy.
Note that in all the datasets, \vesta successfully finds sparsity points at the elbows of the RE curves, resulting in reasonable sparsity and RE trade-offs.
% In all three datasets, accuracy curve had a left flipped `L' shape.
Similar trend was observed for the sparsity and the test RE (see the supplementary material\cite{supp-material}).

\subsection{Data Scalability}\label{sec:scalablility}
\begin{figure}[ht!]
    \centering
    \includegraphics[width=\textwidth]{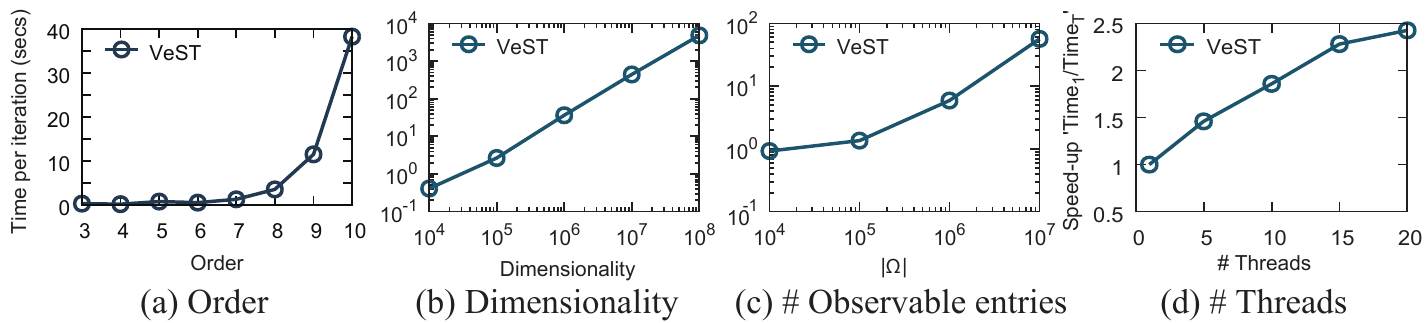}
    \caption{Scalability of \vest. Data scalability of \vest varying
    (a) tensor order, (b) tensor dimensionality, (c) number of observable entries, and (d) number of threads.
    }
    \label{fig:scalability}
\end{figure}
We evaluated the scalability of \vest by
generating synthetic tensors varying the order, the dimensionality, and the number of observable entries, and measuring the running time (see Figure~\ref{fig:scalability}).
%The entry values of the synthetic tensors are randomly selected real numbers between 0 and 1.
For convenience, the dimensionality of each mode in the input tensor, as well as the dimensionality of core tensor, were set equal, i.e., $I_1 = I_2 = ... = I_N$ and $J_1=J_2=...=J_N=3$, respectively.

\noindent \textbf{Order.}
Data scalability on the order of input tensor is tested on synthetic tensors of varying orders from 3 to 10.
For each input tensor, the dimensionality and the number of observable entries were fixed to $|\Omega|= 10^3$.
Figure~\ref{fig:scalability} (a) shows that \vest scales quadratically with regard to order, as discussed in Section~\ref{method:subsec:parallel}.

\noindent \textbf{Dimensionality.}
Data scalability on the dimensionality was tested on input tensors of varying dimensionalities from $10^3$ to $10^8$ with each mode having equal dimension.
The order was set to three, and $|\Omega|$ was set equal to the dimensionality.
Figure~\ref{fig:scalability} (b) shows that \vest has near-linear scalability in terms of the dimensionality.

\noindent \textbf{Number of Observable Entries.}
Data scalability on the number of observable entries was tested on input tensors by varying the number of observable entries from $10^3$ to $10^7$.
The order was set to three, and the dimensionality was fixed to $10^3$.
Figure~\ref{fig:scalability} (c) shows
that \vest has near-linear scalability in term of the number of observable entries.

\noindent \textbf{Effectiveness of Parallelization.}
We evaluated the parallelization scalability of \vest by increasing the number of threads from 1 to 20 and measuring $Time_1/Time_T$ where $Time_T$ is the running time per iteration using $T$ threads.
Figure~\ref{fig:scalability}(d) shows near-linear scalability of \vest in terms of the number of threads used.

\subsection{Discovery}\label{sec:discovery}
We evaluated interpretability of \vest by investigating the factorization results of MovieLens dataset and visually showing that the sparse results enhance interpretability.
% interpretation of sparse results provided by \vest is easier compared to dense results.
It is difficult to analyze dense results generated by vanilla methods without post-processing.
In contrast to existing methods, we can easily identify interesting factors generated by \vest based on sparsity of each row of a factor matrix.

\begin{figure}[h!]
    \centering
    \includegraphics[width=0.9\textwidth]{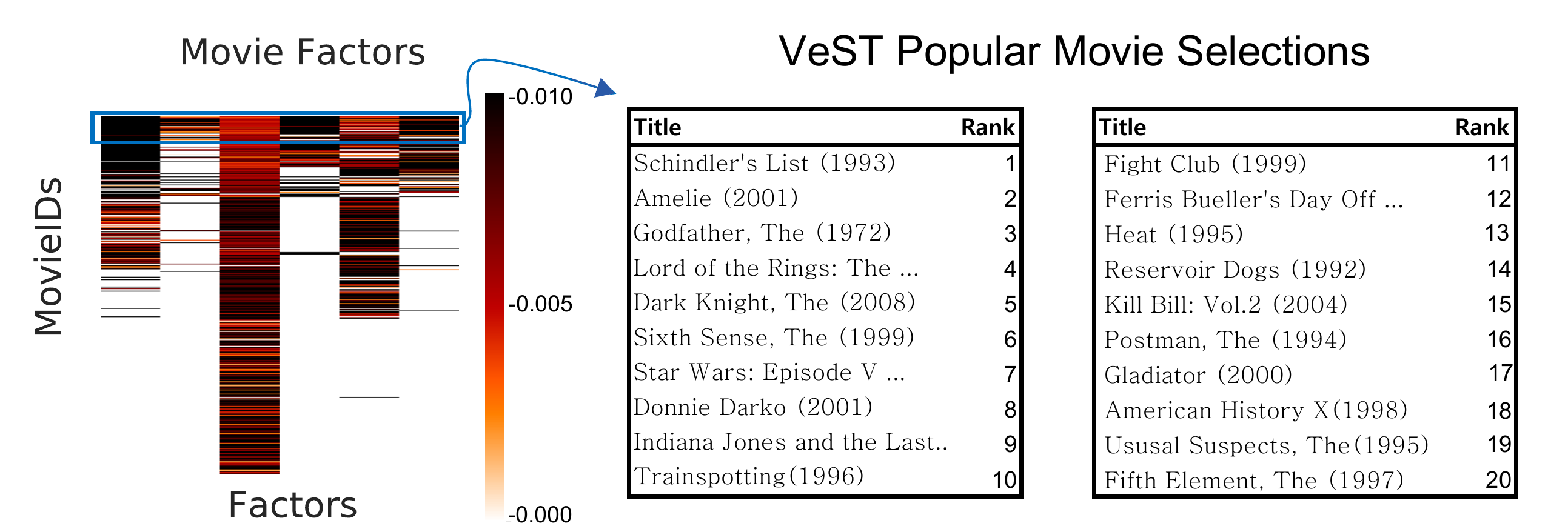}
    \caption{Popular movies discovered by \vest. 14 of 20 movies we find are included in the 100 greatest movies introduced in Empire magazine.
    The remaining six movies (Sixth Sense, Kill Bill, Ferris Bueller's Day Off, Postman, American History X, and Fifth Element) are also famous movies and were rated by various users.
    }
    \label{fig:discovery}
\end{figure}

\textbf{Discovery of Greatest Movies.}
We found that a few rows of the movie-associated factor matrix are fully dense although the goal of \vest is to generate sparse results.
Such rows corresponded to popular movies rated by diverse users.
%
% First, those movies are linked with various ratings of MovieLens data when all factors of movies are non-zero.
% For instance, those movies are rated by users who have different taste, or rated for a long time.
% Also, we interpret the sum of values of all factors as the influence of the movie by considering the value of each factor as a score.
%
Figure~\ref{fig:discovery} shows popular movies which have the largest number of non-zero and the sums of values.
According to Empire magazine~\cite{empireMagazine}, 14 out of 20 movies we found were included in the 100 greatest movies.
The remaining six movies, including `Sixth Sense', `Kill Bill', and `Fifth Element', that were not in the 100 list were also very famous.

%	
%\section{Discovery}
%\label{sec:discovery}
%\input{050discovery}
%
%\section{Related Works}
%\label{sec:related_works}
%\input{060related}
%
\section{Conclusion}
\label{sec:conclusion}
We proposed \vest, a very-sparse Tucker factorization method for sparse and partially observable tensors.
By deriving the element-wise partial differential equations, determining the importance of elements by responsibilities, and parallel distribution of computational work, \vest successfully offers very sparse and accurate results that are applicable for large partially observable tensors.
\vest generates at least 2.2 times more sparse results compared to other methods for partially observable tensors and at least 2.8 times accurate results compared to other sparse factorization methods.
\vest also shows near linear scalability regarding tensor dimensionality, number of observable entries, and number of threads.
Thanks to the increased sparsity that leads to improved interpretability by \vest, we were able to discover interesting patterns related to the greatest movies in the factor matrix of a real-world movie rating tensor data.
Future works include better initialization for Tucker factorization, integration of prior knowledge, and effective visualization of tensor results.

\bibliographystyle{splncs04}
\bibliography{BIB/jeong}

% % ----- APPENDIX ----
% \clearpage
% \newpage
% \setcounter{figure}{0}
% \setcounter{table}{0}
% \setcounter{section}{0}
% \setcounter{equation}{0}
% \appendix
% \input{080appendix}

\end{document}